\documentclass[12pt]{article}
\usepackage{bbm} 
\usepackage[centertags]{amsmath}
\usepackage{amsfonts}
\usepackage{amssymb}
\usepackage{amsthm}
\usepackage{amsmath,mathrsfs}
\usepackage{amsmath,amscd}
\usepackage{mathtools}
\usepackage[matrix,arrow,curve,cmtip]{xy}
\title{\textbf{Introduction to the Category of Derived Motivic Spectra}}
\author{Renaud Gauthier \footnote{2020 Math. Subj. Class: 14C15, 18N60, 55P42, 18N55, 14F42. Keywords: Motives, Segal topoi, derived stacks, spectra.} \\ \\}
\theoremstyle{definition}

\newcommand{\beq}{\begin{equation}}
\newcommand{\eeq}{\end{equation}}

\newcommand{\hrarr}{\hookrightarrow}
\newcommand{\rarr}{\rightarrow}

\newcommand{\rlarr}{\rightleftarrows}
\newcommand{\larr}{\leftarrow}

\newcommand{\eset}{\emptyset}
\newcommand{\Ob}{\text{Ob\,}}

\newcommand{\xrarr}{\xrightarrow}

\newcommand{\rrarr}{\rightrightarrows}

\newcommand{\cA}{\mathcal{A}}

\newcommand{\cC}{\mathcal{C}}
\newcommand{\cCop}{\cC^{\op}}
\newcommand{\cD}{\mathcal{D}}
\newcommand{\cE}{\mathcal{E}}
\newcommand{\cF}{\mathcal{F}}
\newcommand{\cG}{\mathcal{G}}
\newcommand{\cH}{\mathcal{H}}

\newcommand{\cI}{\mathcal{I}}

\newcommand{\cM}{\mathcal{M}}

\newcommand{\cS}{\mathcal{S}}
\newcommand{\cT}{\mathcal{T}}

\newcommand{\cW}{\mathcal{W}}
\newcommand{\cX}{\mathcal{X}}

\newcommand{\bA}{\mathbb{A}}

\newcommand{\bL}{\mathbb{L}}
\newcommand{\bP}{\mathbb{P}}
\newcommand{\bQ}{\mathbb{Q}}
\newcommand{\bR}{\mathbb{R}}

\newcommand{\bZ}{\mathbb{Z}}

\newcommand{\Cat}{\text{Cat}}

\newcommand{\diag}{\text{diag}}

\newcommand{\Fun}{\text{Fun}}

\newcommand{\Hom}{\text{Hom}}
\newcommand{\Ho}{\text{Ho}\,}

\newcommand{\Map}{\text{Map}}
\newcommand{\map}{\text{map}}

\newcommand{\op}{\text{op}}

\newcommand{\Set}{\text{Set}}
\newcommand{\Spec}{\text{Spec\,}}
\newcommand{\Top}{\text{Top}}
\newcommand{\uHom}{\underline{\Hom}}

\newcommand{\Ab}{\text{Ab}}

\newcommand{\CAlg}{\text{CAlg}}

\newcommand{\Catinf}{\Cat_{\infty}}

\newcommand{\CatD}{\Cat_{\Delta}}

\newcommand{\CommC}{\text{Comm}(\cC)}

\newcommand{\dkAff}{\text{d}k\text{-Aff}}
\newcommand{\dSt}{\text{dSt}}

\newcommand{\DkAff}{\text{D}^- k \text{-Aff}}

\newcommand{\kAlg}{k\text{-Alg}}

\newcommand{\Mhat}{M^{\wedge}}
\newcommand{\Mtildetau}{M^{\sim, \tau}}

\newcommand{\oT}{\otimes}

\newcommand{\Ruh}{\mathbb{R}\uh}

\newcommand{\RuHom}{\bR \uHom}
\newcommand{\RuSpec}{\bR \uSpec}

\newcommand{\sPr}{\text{sPr}}

\newcommand{\Sh}{\text{Sh}}

\newcommand{\SetD}{\Set_{\Delta}}

\newcommand{\Sp}{\text{Sp}}

\newcommand{\skMod}{\text{s}k\text{-Mod}}

\newcommand{\skAlg}{\text{s}k\text{-Alg}}

\newcommand{\Sch}{\text{Sch}}

\newcommand{\uh}{\underline{h}}

\newcommand{\uSpec}{\underline{\Spec}}

\newcommand{\dSch}{\text{dSch}}
\newcommand{\DkAfftet}{\DkAff^{\;\sim \, , \, \acute{e}t.}}
\newcommand{\Shet}{\text{Sh}_{\acute{e}t}}
\newcommand{\Schk}{\Sch_k}
\newcommand{\SchkX}{(\Schk)_{/X}}
\newcommand{\dSchk}{\dSch_k}

\newcommand{\dStkF}{(\dSt_k)_{/F}}
\newcommand{\DM}{\text{DM}}
\newcommand{\DMB}{\DM_B}
\newcommand{\DAone}{\text{D}_{\bA^1}}
\newcommand{\ShNis}{\Sh_{Nis}}
\newcommand{\uDM}{\underline{\DM}}
\newcommand{\Sm}{\text{Sm}}
\newcommand{\SmX}{\Sm_{/X}}
\newcommand{\WOmAone}{\cW_{\Omega, \bA^1}}
\newcommand{\oTL}{\oT^{\bL}}
\newcommand{\EMod}{\cE \text{-Mod}}
\newcommand{\DMtQ}{\DM_{t,\bQ}}
\newcommand{\uDMtQ}{\uDM_{t,\bQ}}
\newcommand{\Sft}{\cS^{ft}}
\newcommand{\DMQ}{\DM_{\bQ}}
\newcommand{\Shtr}{\Sh^{tr}}
\newcommand{\SmcorQX}{\Sm^{cor}_{\bQ,X}}
\newcommand{\DAoneet}{\text{D}_{\bA^1,\acute{e}t}}
\newcommand{\HB}{\text{H}_B}
\newcommand{\Shtau}{\Sh_{\tau}}
\newcommand{\cCX}{\cC_{/X}}
\newcommand{\SH}{\cS \cH}
\newcommand{\SHT}{\SH_{\cT}}
\newcommand{\SHPone}{\SH_{\bP^1}}
\newcommand{\SmS}{\Sm_{/S}}
\newcommand{\Spc}{\text{Spc}}
\newcommand{\WhatdkAff}{\widehat{\dkAff}}
\newcommand{\uAonek}{\underline{\bA^1_k}}

\newcommand{\RHomdStkFopTop}{\RuHom(\dStkF^{\op},\Top)}
\newcommand{\Lmot}{\text{L}_{mot}}
\newcommand{\LAone}{\text{L}_{\bA^1}}

\newcommand{\LNis}{\text{L}_{Nis}}
\newcommand{\LuAonekNis}{\text{L}_{\uAonek \cup Nis}}
\newcommand{\dStkFop}{\dStkF^{\op}}
\newcommand{\PoneX}{\bP^1_X}
\newcommand{\PoneiX}{\bP^1_{iX}}
\newcommand{\hocolim}{\text{hocolim}}
\newcommand{\holim}{\text{holim}}
\newcommand{\Mot}{\text{Mot}}
\newcommand{\SeMot}{\text{SeMot}}
\newcommand{\MotSp}{\text{MotSp}}

\begin{document}
\maketitle
\begin{abstract}
	We formalize an abstraction of Grothendieck's philosophy of motives and construct a category of derived motivic spectra in the Segal category $\RHomdStkFopTop$ ($\dSt_k$ the Segal category of derived stacks on $\skAlg$, $\Top = L\SetD$ the Segal category of simplicial sets), thereby providing a starting point for a construction of a stable motivic homotopy category in the Segal setting.
\end{abstract}

\newpage

\section{Introduction}
Grothendieck observed that in algebraic geometry one has a plethora of seemingly different cohomology theories that sometimes yield the same groups, which led him to think that perhaps there ought to be a common underlying theory to all these cohomology theories, something he dubbed motives (\cite{Gr}), an object intermediate between a space and a cohomology. Formalizing this idea, one can construct derived categories of sheaves in an attempt to construct categories of motives (\cite{CD}), which if well-defined enough have objects within them that represent cohomologies, and if in addition those categories are endowed with a 6 functors formalism, which is in a sense a shadow of cohomological properties, one should regard those categories as categories of motives. Already at this point one has shifted focus from cohomologies to derived categories of motives that predate cohomologies, something we will exploit to motivate our approach. Now Khan used the fact that many problems in classical algebraic geometry can be simplified in derived algebraic geometry (most notably for us intersection theory, which involves algebraic cycles, themselves needed to define Weil cohomologies, which are intimately tied to a notion of motives), to warrant developing a theory of motives in that context, which he did in \cite{K}. We will follow the same point of view and use Segal categories for that purpose.\\

Most importantly, what we do presently is expose an alternative way of thinking about motives and how they come about. If classically motives $hX$ associated with a space $X$  appear as in the following diagram:
\beq
\xymatrix{
	X \ar[r] & hX \ar[dl] \\
	H^*(X)
} \nonumber
\eeq
where $H^*(X)$ denotes a Weil cohomology, we suggest what initially seems to be a functorial formalism as in:
\beq
\xymatrix{
	X \ar[d] \ar[r] & hX \ar[dl] \ar@{.>}[d] \\
	H^*(X) & H^*(hX) \ar@{.>}[l] 
} \nonumber
\eeq
where $H^*(hX)$ would capture the cohomological information from a motive. Since, as we argued above, our attention can shift from cohomologies to studying derived categories, we seek diagrams of the following form instead:
\beq
\xymatrix{
	X \ar[d] \ar[r] & hX \ar[dl] \ar@{.>}[d] \\
	\text{D}(X) & \text{D}(hX) \ar@{.>}[l] 
} \label{dercatCD}
\eeq
We do have something formally similar to this picture where derived categories of sheaves are being used, in which setting the bottom map above can be regarded as a specialization functor, thereby prompting us to view the top map as a generalization. Indeed, from \cite{CD} for instance we have:
\beq
\xymatrix{
	\Sch^{ft,sm}_{/X} \ar[d] \ar[r] & \Sch^{ft}_{/X} \ar[d] \\
	\DAone \Shet (\Sch^{ft,sm}_{/X},\bQ)  & \DAone \Sh_t (\Sch^{ft}_{/X},\bQ) \ar[l]
} \nonumber
\eeq
where $t = h$ or $qfh$. Coming back to \eqref{dercatCD}, we regard the derived categories of sheaves therein as cohomological manifestations of schemes $X$ and corresponding motives $hX$. To be fair, those are functor categories having slice categories $\cC_{/X}$ and $\cC'_{/hX}$ for base categories. This prompts us to regard the philosophy of motives as a multi-layered process, the first step consisting in finding a large ambient category $\cC'$ an object $X$ of $\cC$ can be embedded into and putting that object in relation with other objects of that category by  means of a categorical perception such as a slice category $\cC'_{/hX}$, the idea being that the object $X$ would be fully displayed in this setting. The second step consists in turning slice categories into sites. We regard putting a Grothendieck topology on categories as providing a \textbf{filter} on such categories, the idea being that appropriate filters would put objects of such categories in a favorable light. The next step consists in studying manifestations of such sites. The proper way to do this coherently of course is to use sheaves, hence the use of derived categories of sheaves. Finally to highlight a particular characteristic of a geometric object $X$ one then has to select an appropriate modification of such derived categories to highlight such a characteristic, something we regard as a manifestation thereof. We refer to this last step as a \textbf{reduction}.\\ 

\newpage

Thus we are considering multi-layered diagrams of categories, each step of which being regarded as a filtering process, an example being the following diagram:
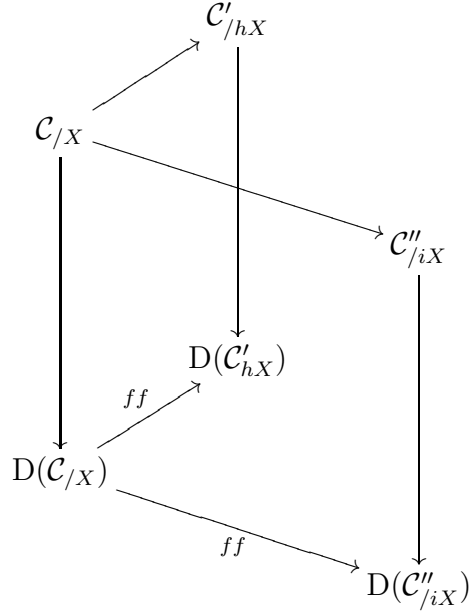
\beq
\xymatrix{
	& \cC'_{/hX} \ar[ddd] \\
	\cC_{/X} \ar[ddd] \ar[ur] \ar[rrd]  \\
	&& \cC''_{/i X} \ar[ddd] \\
	& \text{D}(\cC'_{hX}) \\
	\text{D}(\cC_{/X}) \ar[ur]^{ff}  \ar[rrd]_{ff}\\
	&& \text{D}(\cC''_{/iX})
} \nonumber
\eeq
each vertical map connecting one layer with the next. Such diagrams are functorial in $X$.\\

For the sake of constructing such multi-layered diagrams, in their first stages we regard embeddings as a fundamental initial step towards highlighting the various facets of a geometric object $X$. Given an embedding $h: \cC \hrarr \cC'$ with $X \in \cC$ and $hX \in \cC'$, the categorical perceptions figuring in the induced map $\cC_{/X} \hrarr \cC'_{/hX}$ allow one to see how much of the original object $X$ is perceived at each stage. In the language of derived categories of \eqref{dercatCD}, if for any other embedding $i: \cC' \hrarr \cC''$ no information is gained by studying the embedding of $hX$ into $\cC''$, that is $\text{D}'(hX) \cong \text{D}''(ihX)$, where $\text{D}'(hX)$ is a derived category of sheaves on $\cC'_{/hX}$, $\text{D}''(ihX)$ one on $\cC''_{/ihX}$, and this for any $\text{D}''$, then we say $\text{D}'(hX)$ fully exhibits the nature of $X$ via $hX$, or that $hX$ is the \textbf{complete} object associated with $X$. To reach this point already one has to use the appropriate filters and manifestations on such categories $\cC'$ and $\cC''$, each such modification adding a corresponding layer to our construction. If one wants to highlight a particular characteristic $T$ of a geometric object $X$, if $\mathfrak{R}_T$ designates its corresponding reduction on derived categories of sheaves for instance, if after reduction we have a similar statement as in $\mathfrak{R}'_T \text{D}'(hX) \cong \mathfrak{R}''_T \text{D}''(ihX)$ for any $\mathfrak{R}''_T$, then one says $\mathfrak{R}'_T \text{D}'(hX)$ exhibits the $T$-character of $X$ via $hX$, it captures its $T$-ness. By convention, if $T = \eset$, $\mathfrak{R}_{\eset} = id$, the object $X$ has been fully displayed in $hX$, which consequently becomes the complete object associated with $X$. If one considers those reductions $\mathfrak{R}_{mot}$ that allow us to obtain motivic categories, namely those that turn derived categories of sheaves of the form $\DAone \Sh_{\tau}(\cC_{/X}, \bQ)$ as in \eqref{dercatCD} for example into categories for which a 6 functors formalism holds, from our perspective the resulting categories $\mathfrak{R}_{mot} \DAone \Sh_{\tau}(\cC_{/X}, \bQ)$ exhibit the motivic character of $X$. This is a formal heuristic take on the philosophy of motives that we find helfpul to put things in perspective. Nevertheless for the sake of constructing a Segal category of derived motivic spectra, which is the main goal of this paper, we will base our work on the notion of motives as they are usually understood. Our heuristics on motives will be used to smoothly introduce the Segal category derived motivic spectra are issued from.\\

We start with a scheme $X \in \Schk$, that we first embed in derived schemes, $iX \in \dSchk$, and then derived stacks using the model Yoneda embedding, $F = h_{iX} \in \dSt_k$ (see \cite{HAGII}). We then focus on the slice category $\dStkF$. If most derived categories of sheaves in the traditional study of motives are of the form $\text{D} \Shtau(\cCX,\cA)$, we argue the correct generalization of such an object in our context is constructed from a Segal category $\RHomdStkFopTop$. We develop a formalism of derived motivic spectra of objects in this category, the equivalent in our setting of motivic spectra in \cite{K}. We define derived motivic spectra to be spectra that are $\bA^1$-invariant (or rather its generalization) and satisfy Nisnevich excision instead of Nisnevich descent. One thing we do not do however is develop a 6 functors formalism, and we do not work with pointed spaces for simplicity.\\

About notations, a hom set is denoted with the letter $\Hom$, $\uHom$ is used for internal homs, $\Map$ for simplicial mapping spaces, and $\map$ for a homotopy function complex. $\CatD$ will denote the category of simplicial categories and $\SetD$ the category of simplicial sets. We refer to $\Hom(-,X)$ as a functorial perception of $X$, $\cC_{/X}$ as a categorical perception, $\Hom(X,-)$ as a functorial manifestation of $X$ and $\cC_{X/}$ as a categorical manifestation.\\

\section{From embeddings to comparison theorems}

\subsection{Concrete realization vs. abstraction}
Working with motives (\cite{Gr}, \cite{S}, \cite{A}) means in particular working with Weil cohomologies (\cite{J}), since those factor through motives, and Weil cohomologies necessitate that we consider algebraic cycles among other things. This forces the standard conjectures upon us, which impose very concrete demands on our cohomology theories, which consequently appear as realizations of the overarching theory of motives they are issued from (\cite{M2}). It follows that developing a theory of motives with Weil cohomologies in mind provides a facet only of a theory of motives. Adopting such a point of view subsumes motives are fundamental objects that take preeminence over cohomology theories. Thus what we are aiming for is a precursor to the notion of classical motive from which cohomologies can be derived.\\

To illustrate the need to seek an abstract theory of motives, recall that the concept of motive $hX$ associated with a geometric object $X$ was introduced since various cohomologies $H^*_{(i)}(X)$ could produce the same results, thereby exhibiting those cohomologies as representatives of a same class, the motive $hX$. We have the following picture:
\beq
\xymatrix{
	X \ar[r] & hX \ar[dl] \\
	H^*_i(X)
} \label{triangle}
\eeq
Given that a motive is an object intermediate between a space and a cohomology theory, suppose for a moment we had a way of dissociating the cohomological information from a motive by means of a functor $H^*$ in such a manner that \eqref{triangle} could be completed into a square diagram as follows:
\beq
\xymatrix{
	X \ar[d] \ar[r] & hX \ar[dl] \ar@{.>}[d] \\
	H^*_i(X) & H^*(hX) \ar@{.>}[l]
} \nonumber
\eeq
where each diagonal map factors through specialization functors $H^*_i(X) \larr H^*(hX)$, thereby presenting $X \rarr hX$ as a generalization. That we have such a formalism is not as far fetched as it seems since we will see later that cohomological spectra are objects of derived categories of sheaves (or stacks) on certain slice categories and that the above statement about cohomological objects should read $\text{D} \Sh_{\tau}(\cC_{/X},\cA) \larr \text{D} \Sh_{\tau'}(\cC'_{/hX}, \cA)$, a pre-comparison functor of some sort, which follows from having a map $X \rarr hX$. Making such a statement already shows a clear departure from standard motivic homotopy theory whereby $\text{D}\Shtau(\cC_{/X},\cA)$ is a category of motives, while for us the motivic character of a geometric object $X \in \cC$ is apparent in its embedding $hX$ in another category $\cC'$ as discussed in the introduction, the various derived categories on the corresponding slice categories being different ways to make manifest the motivic character of such geometric objects. The rest of this paper has as one of its aims to present this alternate formalism that provides an impetus for a theory of motives, in the hope it will spur new research.\\

\subsection{Embeddings}
If a motive is constructed starting from a variety, or a scheme, we argue the ``motivic heart" of a geometric object can be revealed once the said object is first embedded in a larger category, and it is put in contrast with other objects of this category. Thus as a first step into abstraction proper, we are looking at a series of embeddings:
\beq
\text{Var}_k \hrarr \Sch_k \xrarr{i} \dSchk \xrarr{h} \dSt_k \label{emb1}
\eeq
for $k$ a commutative ring, from varieties to schemes, if one is starting from varieties, then derived schemes, followed by derived stacks. To be precise, if $k$ is an algebraically closed field, we have a fully faithful functor $\text{Var}_k \hrarr \Sch_k$. From \cite{DAG}, we have a fully faithful functor from the category of $k$-schemes $\Sch_k$ to the Segal category $\dSchk$ of derived schemes over $k$, which itself is embedded in the Segal category of derived stacks $\dSt_k$. To be more precise, we have a Yoneda embedding $h: \dSchk \rarr \dSt_k$ that sends a derived scheme $iX$, $X \in \Schk$, to $h_{iX} = \Map_{\dSchk}(-,iX)$.\\

Regarding the latter embedding, recall from $\cite{HAGII}$ that if $M$ is a model category, $\Gamma_*: M \rarr M^{\Delta}$ a cofibrant resolution functor on $M$, if we let $\uh_x(y) = \Hom_M(\Gamma_*y,x) = \Map_M(y,x)$, the functor $\uh: M \rarr \Mhat$ (where $\Mhat$ is the left Bousfield localization of $\sPr(M)$ with respect to weak equivalences in $M$) can be right derived into $\Ruh$. Further we have $h_x \cong \Ruh_x$ in $\Ho(\Mhat)$. Now to argue $h_x$ is a stack it suffices to show $\Ruh_x$ is one as well. Recall that a model pre-topology $\tau$ on $M$ is said to be subcanonical if $\Ruh_x$ is a stack for all $x \in M$. Formally this means $\Ruh: \Ho(M) \rarr \Ho(\Mhat)$ factors through the category of stacks $\Ho(\Mtildetau)$. Since this is what we need, we require that $\tau$ be subcanonical. Specifically, here we work with $M = \DkAff = \skAlg^{\op}$. Now that $\tau$ is subcanonical in general follows from Assumption 1.3.2.2 of \cite{HAGII}, which holds in the special case of having $\tau =$ \'et on $\DkAff$. Thus $\Ruh: \Ho(\DkAff) \rarr \Ho(\DkAfftet)$ is fully faithful, which is isomorphic to $h$, thus when we say $h: \dSchk \rarr \dSt_k = L \DkAfftet$ after simplicial localization of model categories is fully faithful, we really mean this in a homotopical sense.\\ 

For completeness' sake, recall how $\dSt_k$ is constructed; one starts, as recalled in \cite{HAGII} and \cite{HidSt}, with a symmetric monoidal model category $\cC$, and in derived algebraic geometry we are really thinking of $\cC = \skMod$. Then $\CommC = \skAlg$. Taking the Dwyer-Kan simplicial localization of this (\cite{DK}) produces a simplicial category $L(\skAlg)$, hence a Segal category (\cite{SeT}). Its opposite Segal category is denoted $\dkAff = L(\skAlg)^{\op}$. Denoting $L\SetD$ by $\Top$, we consider the Segal category of prestacks $\WhatdkAff = \RuHom(\dkAff^{\op}, \Top)$, and $\dSt_k \subset \WhatdkAff$ is the sub-Segal category of derived stacks for the \'etale topology, which can also be obtained as $\dSt_k = L(\DkAfftet)$ where $\DkAfftet$ is the model category of derived stacks on $\skAlg$, as used in the preceding paragraph.\\

\eqref{emb1} provides a chain of embeddings (that can possibly be further extended), and we seek the most general category a geometric object $X$ can be embedded in, the idea being that if it appears difficult to find an optimal category $\cC'$ such that any other embedding $\cC' \hrarr \cC''$ would not provide any gain in information about the object $X$ as discussed in the introduction, assuming the $T$-ness of $X$ is already present in $hX \in \cC'$, any further embedding thereafter would preserve the $T$-character of $X$, thus working with the most general category one can think of should suffice. Presently we will settle on $\dSt_k$ for the sake of defining motives, which is sufficient since the category of smooth schemes already provides a base for constructing categories of motives. Putting $X$ in contrast with other objects of this larger category $\dSt_k$ really means looking at a functorial perception $\Hom(-,X)$, or a categorical perception such as a slice category $(\dSt_k)_{/X}$.\\

To fix ideas, using \eqref{emb1} above, for $X \in \Schk$, an embedding of slice categories is of the form:
\beq
\SchkX \hrarr (\dSchk)_{/iX} \hrarr (\dSt_k)_{/h_{iX}} \label{emb2}
\eeq
The motivation for working with slice categories follows from the relative point of view philosophy of Grothendieck according to which the various aspects of an object can be made manifest once it is put in relation with other objects of the category in which it finds itself, something we refer to as \textbf{perception}. In particular \eqref{emb2} amounts to working with categorical perceptions. Next one puts topologies on such slice categories, something we regard as filters. Already at this point we are starting to have layered diagrams, with slice categories as base layer, followed by a layer of filtered categories, or sites. Then one studies sheaves on such sites, since they provide coherent manifestations, and this produces an additional layer. Another way of looking at this is that categories of sheaves formally correspond to an extrinsic comparison of a given base category with various target categories, something we will discuss in the next section. Finally one studies appropriate functors on the resulting topos to study certain aspects of $X$, something we referred to as reductions in the introduction. If one desires to focus on the motivic aspect of $X$, one works with reductions producing cohomological functors, and by this we mean functors that satisfy a 6 functors formalism, which in a sense extract from $X$ its cohomological character. This is in essence the program we will follow.

\subsection{Manifestations/ extrinsinc representations}
Following Grothendieck's philosophy of working with relative objects, we have introduced categorical perceptions $\SchkX$ and $\dStkF$, with $F = h_{iX}$. Categorical perceptions and their relative embeddings are primordial in our philosophy insofar as once appropriately chosen they provide a foundation for determining at which point a given geometric object has its $T$-ness made manifest. As briefly mentioned at the end of the preceding section, the next step consists in putting a topology on those slice categories, a filtering process of some sort, and we assume that has been done at this point.\\

\newpage

Still following Grothendieck's relative point of view philosophy, analyzing those sites should involve other categories. One way to implement this is by studying their extrinsinc representations, and this we can do functorially using sheaves for instance. Those functorial manifestations when appropriately modified via reductions can exhibit one aspect or another of geometric objects, and it is the resulting reduced manifestations that tell us whether we have reached an optimal ambient category $\cC'$ in a tower of embeddings, or whether we have to push further in such a tower.\\

This is actually very general: to highlight the character of type $T$ of a geometric object $X \in \cC$, one seeks ambient categories $\cC'$ with corresponding embeddings $i: \cC \hrarr \cC'$ such that for appropriately chosen functors $F_T$ one says the $T$-character of $X$ is made manifest in $iX \in \cC'$ if for any other category $\cC''$ and any other embedding $j: \cC' \hrarr \cC''$, for any other $T$-functor $G_T$, $F_T(\cC'_{/iX}) \cong G_T(\cC''_{/jiX})$, at which point we would also say $F_T(\cC'_{/iX})$ displays the type $T$ of $X$. If $T$ is being a motive and $F_T$ corresponds to a derived category of sheaves on $\cC'_{/iX}$ satisfying a 6 functors formalism, in this setting $F_T(\cC'_{/iX})$ is a category of motives in the traditional sense, which in our formalism exhibits the motivic character of $X$, or in other words it captures the motivic aspect of $iX$.\\

Now to mirror the fact that in \cite{CD} we work with sheaves over slice categories, which in themselves are manifestations of such base categories, we now seek in our present context some manifestations of $\dStkF$, and those will be objects of a Segal category $\RHomdStkFopTop$. This will be the analog for us in the Segal setting of the categories of sheaves used in \cite{CD}. Note that we will not put topologies on $\dStkFop$ since this would put a filter on it, and we desire to work in wide generality.\\

\subsection{A study of cohomology - deconstruction}
Recall that one can associate to compact manifolds various cohomology theories that yield the same groups, supposing they are based on the Eilenberg-Steenrod axioms. This presents those cohomology theories as realizations of those axioms, which therefore appear as more fundamental, more motivic in a sense. This also shows that we should distance ourselves from cohomology theories proper since they appear secondary, and aim instead for precursors to cohomology theories, chief among them being motives. To further illustrate the necessity of a change of focus from cohomology to more abstract theories, observe that there are situations where cohomology is instrumentalized. An example of such a situation is sheaf cohomology. As recounted in \cite{M} for example, recall that if $\Shet(X)$ denotes the category of \'etale sheaves on a scheme $X$, the functor $\Shet(X) \rarr \Ab$, $\cF \mapsto \Gamma(X,\cF) = \cF(X)$ is in a sense a first order manifestation of $X$ into $\Ab$, relative to $\cF$ of course. Since the global sections functor $\Gamma(X,-)$ is left exact, one can find its right derived functor. To do this,  since $\Shet(X)$ has enough injectives, one can find an injective resolution $0 \rarr \cF \rarr \cI^0 \rarr \cI^1 \rarr \cdots $ of $\cF$. Because the complex $\Gamma(X,\cI^0) \rarr \Gamma(X,\cI^1) \rarr \cdots $ is in general not exact, one can find its cohomology, which is none other than the \'etale cohomology of $X$. Note here that studying the cohomology of $\Gamma(X,\cI^{\bullet})$ complements its non-exactness, it restores a semblance of coherence. What transpires is that cohomology is necessary to study the manifestation of $X \in \Schk$ into $\Ab$ via $\cF \in \Shet(X)$, for the sake of having a complete picture of $X$, but does not take center stage otherwise. Thus if one desires to develop a theory of motives, regarded as fundamental objects, it seems fitting to not directly use cohomology, which seems at times to be secondary, but to base a theory of motives on more abstract theories instead.\\

\subsection{Shift from cohomology to studying derived categories}
As we pointed out in the preceding subsection, cohomology from our perspective is useful to study geometric objects. However to formalize a theory of cohomologies, one ought to aim for theories that predate cohomologies. Motives is one such theory. The way motives were originally introduced, they were presented as an object in between a space and a cohomological object, in that any good cohomology would factor as in:
\beq
\xymatrix{
	X \ar[dr] \ar[r] & hX \ar@{.>}[d] \\
	& H^*(X)
} \nonumber
\eeq
That motives $hX$ appear thus necessitates that they carry some cohomological information, and for this it is typically required that motives be objects of a category satisfying a 6 functors formalism. This is not something we focus on, we will not ask that our categories of coefficients, as they will later be called, satisfy a 6 functors formalism. We will explain later why that is the case. This simplifies our task greatly and to develop a 6 functors formalism is still possible using our categories of derived motivic spectra as a starting point. \\

Categories of coefficients are derived categories that predate cohomology theories. In \cite{CD} in particular, many of those in the context of motives are derived categories of sheaves, and some of them have cohomology spectra as objects that represent cohomology theories. This is our desired abstraction. In this manner one can eschew cohomology theories in favor of studying such derived categories, and comparisons between them. This is typically what is done since those derived categories are categories of motives. In the present paper we broaden the picture a bit by seeing those derived categories as manifestations of geometric objects. As such, they are meant to exhibit the motivic character of such objects, as recounted in the introduction.\\

\subsection{Categories of motives}
We use the presentation given in \cite{CD} of the various categories of motives of interest to us. Central for us is the derived category of Beilinson motives $\DMB(X)$, the Verdier quotient of $\DAone(\ShNis(\SmX,\bQ))$ by a subcategory of acyclic objects. It is worth pointing out again that if this is generally accepted as providing a category of motives, we propose an alternative way of seeing the study of motives as equivalent to first placing a given geometric object $X$ in the most general ambient category it can be seen as an object of, and seeing it relative to other objects in this category. Thus we are looking at slice categories $\cC_{/X}$. Those can be analyzed by using what we referred to as filters, namely Grothendieck topologies that turn them into sites. At this point one projects such sites into other categories using sheaves, thereby providing us with topos. The various facets, or characteristics of type $T$, of the geometric object $X$ we started with can be studied using $T$-specific manifestations called reductions, that hone in on the $T$-character of geometric objects. From that perspective we can regard $\DMB(X)$ as one reduction that provides a manifestation of the motivic character of $X \in \Sm$ by using Nisnevich sheaves with rational coefficients. Observe that this in no way changes all the work done on motives so far, it just provides another approach to the philosophy of motives; if $DM_B(X)$ is a category of motives in the traditional sense, we merely suggest seeing it as exhibiting the motivic character of $X$.\\

To come back to $\DMB(X)$, by $\Sm$ we mean the class of smooth morphism of finite type in a category of schemes over a scheme $S$ fixed at the onset. By $\DAone(\cA)$ for some category $\cA$ we mean $\DAone(\cA) = \text{D}(\Sp(\cA))[\WOmAone^{-1}]$, where $\text{D}(\cC)$ is the usual derived category of some category $\cC$, $\Sp(\cA)$ is a category of spectra in $\cA$, and here the localization is with respect to a class $\WOmAone = \cW_{\Omega} \cup \cW_{\bA^1}$, the collection of two families of morphisms of Tate spectra. The reader is invited to read \cite{CD} for more details, our very cursory presentation is just meant to present the main actors and put things in perspective, we shall not deal with details beyond this level of presentation in the rest of this work. For us what is important to keep in mind is that the derived category of Beilinson motives is a derived category of rational, Nisnevich sheaves on the smooth site of a scheme $X$.\\

We also have the category of Voevodsky's motives $\DMtQ(X)$ where the topology $t = h$ or $qfh$, defined as the localizing subcategory of $\uDMtQ(X)$ spanned by objects of the form $\Sigma^{\infty} \bQ^t_X(S)(n)$ for $S \in \SmX$, $n \leq 0$, where $\bQ^t_X(S)$ is the free $\bQ$-linear $t$-sheaf represented by $S$, $\uDMtQ(X) = \DAone(\Sh_t(\Sft_{/X},\bQ))$, $\Sft$ a class of morphism of finite type in a category $\cS$ of Noetherian schemes over $X$ satisfying certain conditions (see \cite{CD} for details).\\

Another category of motives is $\DMQ(X)$, defined as a derived category of sheaves with transfer $\DAone(\Shtr(X,\bQ)) = \DAone(\ShNis(\SmcorQX))$.\\

The category of Morel motives $\DAone(X,\bQ)_+$ is a subcategory of $\DAone(X,\bQ) = \DAone(\ShNis(\SmX,\bQ))$ obtained by projection (see \cite{CD} for details).\\

\subsection{Realization functors}
We can start our discussion of the concept of realization functors by the following result, as presented in \cite{CD}: for a scheme $X \in \Sch_{/S}$ with structural morphism $f: X \rarr S$, $E$ a given stable cohomology with associated ring spectrum $\cE$, we have a realization functor $\DMB(X) \rarr \Ho(\cE_X \text{-Mod})$, where $\cE_X = \bL f^* \cE$, a map that to $M$ associates $\cE_X \oTL_X M$.\\

In \cite{CD}, an abstraction of this situation is given as follows. If $\cM$ and $\cM'$ are monoidal model categories satisfying a host of conditions as spelled out in \cite{CD}, $\phi^*: \cM \rlarr \cM': \phi_*$ is a left Quillen adjunction, $\cE'$ is a fibrant resolution of $\mathbbm{1}$ in $\cM'(\Spec k)$, $\cE$ a cofibrant resolution of $\phi_*(\cE')$ in $\cM(\Spec k)$, then by Thm 17.1.5 of \cite{CD}, we have a realization functor $\Ho(\cM) \rarr \Ho(\cM') \cong \Ho(\EMod)$, $M \mapsto \cE \oTL M$. Note that having $\Ho(\cM') \cong \Ho( \EMod)$ means the representation theory of $\cE$ is given by $\cM'$, and if $\cE$ represents a cohomology theory, then $\cM'$ appears as the ambient category within which the cohomology theory $\cE$ is supposed to represent can be developed. For this to happen $\cM$ must be a reference category in a sense, and the fact that it provides ``clean" representations of cohomologies such as $\Ho(\EMod)$ justifies that it be considered as being associated with a theory of motives.\\ 

Having presented this abstract set up, observe that if $H_{B,X}$ is the Beilinson motivic cohomology spectrum, then by Thm 14.2.9 of \cite{CD} the functor $\DMB(X) \rarr \Ho(H_{B,X} \text{-Mod})$ is an equivalence of triangulated monoidal categories, thus we have a trivial realization functor illustrating the fact that $\DMB(X)$ is an optimal category of motives. In this manner realization functors have enabled us to justify that derived categories such as $DM_B(X)$ effectively be categories of motives. This is nothing new of course, but in our formalism this has the added benefit of showing that of the two cohomological manifestations $DM_B(X)$ and $\Ho(\cE_X \text{-Mod})$, the former appears as an optimal functor of $X$ that exhibits its motivic character.\\

\subsection{Comparison functors}
Comparison theorems allow us to link all the motivic categories of interest. One thing worth pointing out before doing so is that if we have such equivalences of categories, this indicates that each of the categories involved are not absolute in a sense, but are merely facets of a yet to be determined ultimate category of motives, each of those categories figuring in those comparison theorems being in a sense representatives of the class corresponding to a universal motivic homotopy theory. This point of view justifies the necessity of seeking a higher theory to extract a category of motives all the categories below being representatives of. For the time being though, observe that all the categories of motives involved in the comparison theorems we will cover come from reductions of various derived categories of sheaves as we will see and exhibit the motivic character of $X$ in our formalism.\\

A first comparison theorem is provided by the following, Thm 16.1.2 of \cite{CD}, which states that $\DMB(X) \cong \DMtQ(X)$, where $t = h$ or $t = qfh$. One thing we will be careful to clearly enunciate is what property $X$ is supposed to satisfy in each of the comparison theorems we will give, for a reason we will expound upon below. Here $X$ is an excellent Noetherian scheme of finite dimension.\\

If now $X$ is excellent and geometrically unibranch in addition to being Noetherian and of finite dimension, we have an equivalence of triangulated monoidal categories $\DMB(X) \xrarr{\simeq} \DMQ(X)$.\\

By Thm 16.2.13 of \cite{CD}, for any Noetherian scheme of finite dimension $X$, we have a canonical equivalence of triangulated monoidal categories $ \DMB(X) \simeq \DAone(X,\bQ)_+$.\\

Finally, for any Noetherian scheme of finite dimension $X$, we also have a canonical equivalence of categories $\DMQ(X) \simeq \DAoneet(X,\bQ) = \DAone(\Shet(\SmX,\bQ))$.\\

It is important to note at this point that if we modify the characteristics of our base scheme $X$, we are shifting our point of view, and consequently the motives, understood in the traditional sense, will be different as well. Thus for illustrative purposes, and to fix ideas, we work with one type of scheme only. Thus for all the above comparison theorems to be used in unison, we fix $X$ to be an excellent, geometrically unibranch Noetherian scheme of finite dimension.\\

\section{Abstraction of a category of motivic spectra}
\subsection{General form of reductions} \label{Genformofred}
All the triangulated monoidal motivic categories used in the comparison theorems above are constructed from derived categories of sheaves on the small site of some scheme $X$, thus we will seek an abstract category of motives that retains the same form, namely a category of sheaves, or stacks, on some slice category.\\

To illustrate the fact that the motivic categories we are dealing with are essentially what we call \textbf{reductions} of derived categories of sheaves, we present them in an itemized fashion:
\begin{itemize}
	\item $\DMB(X) = \DAone(\ShNis(\SmX,\bQ))_{/\{ \HB \text{-acyclic} \}} $
	\item $\DMtQ(X) = \langle  \Sigma^{\infty} \bQ^t_X(S)(n) \rangle \subset \DAone(\Sh_t(\Sft_{/X}, \bQ)) $
	\item $\DMQ(X) = \DAone(\ShNis(\SmcorQX)) $
	\item $\DAone(X,\bQ)_+ \subset \DAone(\ShNis(\SmX,\bQ))$ 
	\item $\DAoneet(X,\bQ) = \DAone(\Shet(\SmX,\bQ))$
\end{itemize}
If the symbol $\mathfrak{R}$ designates reduction, all the above categories can be written in the form $\mathfrak{R} \DAone(\Sh_{\tau}(\cCX,\bQ))$, where $\mathfrak{R}$ above is Verdier quotient, projection, subset, or the identity. The categories $\cCX$ are slice categories over a scheme $X$. The larger categories $\DAone(\Sh_{\tau}(\cCX,\bQ))$ are ambient categories where motives in the traditional sense live.\\

\subsection{Multi-layered diagrams of categories}
Now that we have introduced comparison theorems and reductions, we digress momentarily and present an example of a multi-layered diagram. It initially provides slice categories in one layer, followed by sites, then derived categories of sheaves, and in the last layer the reductions of such derived categories to provide categories of motives. Here is one example, where we have used the fact that $DM_B(X) \simeq \DAone \Shet(\Sm_{/X}, \bQ)$:
\beq
\xymatrix{
	& \Sft_{/X} \ar[dd] \\
	\Sm_{/X} \ar[dd] \ar[rr] \ar[ur] && \SmcorQX \ar[dd] \\
	& (\Sft_{/X},t) \ar[dd] \\
	(\Sm_{/X},\acute{e}t) \ar[dd] && (\SmcorQX, Nis) \ar[dd] \\
	& \DAone \Sh_t(\Sft_{/X},\bQ) \ar[dd] \\
	\DAone \Shet(\Sm_{/X},\bQ) \ar[dr]_{\simeq} \ar[ur]^{ff} \ar[rr]_{\qquad ff} && \DAone \ShNis(\SmcorQX,\bQ) \ar[dl]^{\simeq} \\
	& DM_B(X)
} \nonumber
\eeq
where the bottom diagonal maps are equivalences, and the last vertical map is a reduction, consisting into taking a localizing subcategory, followed by an equivalence. The last layer shows that Beilinson motives highlight the motivic character of $X$. We will be working at the level of unrestricted, unreduced derived categories of sheaves however, the next to last layer, which is not shown here.\\

At this point it is fitting to expound on the notion of optimal manifestation of a geometric object. A statement about reductions of the form $F \Sh_{\tau'}(\cC'_{/hX}, \Lambda') \cong G \Sh_{\tau''}(\cC''_{/ihX}, \Lambda'')$ for any embedding $i: \cC' \hrarr \cC''$, for any topology $\tau''$ on $\cC''_{/ihX}$, for any category $\Lambda''$, for any reduction $G$ that produces a derived category of sheaves that satisfies a 6 functors formalism, means $F \Sh_{\tau'}(\cC'_{/hX}, \Lambda')$ is a category of motives in the traditional sense, which from our perspective is optimal in highlighting the motivic character of $X$ through $hX \in \cC'$, it corresponds to an optimal vertical thread in a multi-layered diagram of categories such as the one above. Now we have an ordering on slice categories given by embeddings. Thus the first slice category $\cC'$ for which we have isomorphisms of reductions is optimal in that regard. $\cC'$ being fixed, we also have an ordering on topologies $\tau'$, a coarser topology being optimal if it leads to isomorphisms of reductions. However $\cC'$ and $\tau'$ being fixed, various categories $\Lambda'$ and reductions $F$ may be difficult to compare, so we see it is at the level of categories of sheaves and reductions thereof that we are starting to have classes of categories of motives, absent any clear identification of an optimal category of motives.\\

\subsection{Abstraction of a category of coefficients}
The derived categories of motives introduced in Section \ref{Genformofred} are what Khan refers to as categories of coefficients. Those follow a 6 functors formalism, which abstractly encodes what makes cohomology theories what they are. Hence the latter should be derivable from such categories of coefficients.\\

Since we are constructing categories that will ultimately be used to produce cohomology theories, we are interested in setting up a background from which one can develop a 6 functors formalism. However, our focus is really on the functorial manifestations of motivic objects in our heuristic picture, not so much on cohomology per se. Because of this we do not see an interest in developing a 6 functors formalism for a category of coefficients, and we leave this task to the motivated reader.\\

For $\cC$ an $\infty$-category with fiber products, $\Catinf^{Pres}$ the $\infty$-category of presentable $\infty$-categories, Khan defines a category of coefficients to be given by a symmetric monoidal functor $D: \cCop \rarr \Catinf^{Pres}$. If $\cC$ is a category of schemes, $X \in \cC$, then $D(X)$ is a presentable $\infty$-category whose objects are referred to as sheaves on $X$. The reductions $\mathfrak{R} \DAone(\Sh_{\tau}(\cC_{/X}, \bQ))$ we have considered before are of the form $D(X)$ and are regarded as categories of coefficients by abuse of language (albeit after passing to homotopy categories). We mention categories of coefficients presently since they provide a clear formalization of the concept of motivic category.\\

\newpage

Since we have embeddings $\Schk \hrarr \dSchk \hrarr \dSt_k$, moving a $k$-scheme $X$ to a derived scheme $iX$, and then to a derived stack $F = h_{iX}$, it is therefore natural to consider a category of coefficients to be based on a Segal category of the form $\RHomdStkFopTop$, with $\Top = L \SetD$, $\dSt_k = L(\DkAfftet)$, $\DkAfftet$ the model category of stacks on $\skAlg$, $L$ stands for the simplicial localization of Dwyer-Kan, thereby producing Segal categories, $\dSt_k$ is the Segal category of derived stacks on $\skAlg$, sub-Segal category of $\WhatdkAff = \RuHom(\dkAff^{\op}, \Top)$ where $\dkAff = L(\skAlg)^{\op}$. The reason for having $\Top$-valued functors will be given in Section \ref{Generalization}. \\

Independently Gallauer defines coefficients systems in \cite{Ga}, functors $F: \SchkX^{\op} \rarr \Cat^{st,\oT}_{\infty}$ satisfying certain conditions, which are closely related to the categories of coefficients of Khan in the sense that presentable coefficients systems correspond to motivic categories of coefficients in the sense of Khan. We use Khan's coverage since for our purposes it is easier to generalize to the Segal category setting. But coefficient systems bridge the gap between Khan's work and that of \cite{CD} in the sense that the category of homotopy functors $\Ho(F)$ for $F$ a coefficient system yields in particular what is referred to as motivic triangulated categories in \cite{CD}. Another reason for alluding to Gallauer's work is that the universal coefficient system, which consequently would be the category of motives we are looking for, is none other than the stable motivic homotopy category. We use this fact as an imprint to seek an initial construction of such a category in our setting, before being made 6 functors formalism ready, and we refer to it as the category of derived motivic spectra. We then turn to \cite{K} to construct such a category by using the definition of the stable motivic homotopy category of Khan as motivation.\\

To be clear, Khan defines the stable motivic homotopy category $\cS \cH(X)$ as the category of motivic $\bP^1_X$-spectra which is then enhanced for the sake of making it support a 6 functors formalism. Since again this is not something we are aiming for, we stop at the category of motivic spectra, and thus we are really only looking at a category of derived motivic spectra in our setting.\\

To further highlight the central role played by the stable homotopy category, which we momentarily denote by $\text{SH}$ as in \cite{Ga}, one can connect this category with any coefficient system $C$ using an essentially unique functor $\rho^*: \text{SH} \rarr C$ such that for $S \in \SchkX$, there is a right adjoint to $\rho^*(S)$: $\rho_*(S): C(S) \rarr \text{SH}(S)$ in such a manner that $\rho_*(S) \mathbbm{1}:= \cE \in \CAlg(\text{SH}(S))$ is a motivic spectrum that represents $C$-cohomology.\\

\section{Construction of a category of derived  motivic spectra}

\subsection{Khan's construction}
In order to motivate our definition of a category of derived motivic spectra, we briefly recall Khan's construction from \cite{K}. He defines the category of motivic spectra as $\SH(S) = \SHPone(S)$, the category of motivic $\bP^1_S$-spectra, for $S$ a scheme. Here $\bP^1_S$ is shorthand for $M_S(\bP^1_S)$ (see full definition in Section \ref{GenofP}) pointed at $\infty$, and $\SHT(S)$ denotes the category of motivic $\cT$-spectra, whose objects are $\cT$-spectra over $S$ that satisfy Nisnevich descent and $\bA^1$-homotopy invariance, and this is defined level-wise. Recall that a $\cT$-spectrum is given by a sequence $\{ F_n \}_{n \geq 0}$ of fibred pointed spaces over $S$, along with isomorphisms $F_n \xrarr{\cong} \Omega_{\cT}(F_{n+1})$ for all $n \geq 0$, compatible with morphisms of spectra. Here $\Omega_{\cT}$ is the $\cT$-loop functor $\uHom_S(\cT, -)$ where $\uHom_S$ is the internal hom of the category of smooth fibred spaces over $S$, which are functors $F: \SmS^{\op} \rarr \Spc$ valued in the $\infty$-category of spaces.\\

To satisfy Nisnevich descent for a fibred space $F$ over $S$ is equivalent to satisfying Nisnevich excision, which means that for any cartesian square of smooth derived $S$-schemes of the form:
\beq
\xymatrix{
	U \times_X V \ar[d] \ar[r] & V \ar[d]^p \\
	U \ar[r]_j &X
} \nonumber
\eeq
with $p$ \'etale, $j$ an open immersion, with a complementary closed immersion $Z \hrarr X$ in such a manner that $p^{-1}Z \xrarr{\cong} Z$, the following commutative square in $\Spc$ is cartesian:
\beq
\xymatrix{
	\Gamma(X,F) \ar[d] \ar[r] & \Gamma(U,F) \ar[d] \\
	\Gamma(V,F) \ar[r] & \Gamma(U \times_X V,F)
} \nonumber
\eeq

\subsection{Generalization} \label{Generalization}
Since we have an embedding $\SchkX \hrarr \dStkF$, the role played for us by fibred spaces over a scheme $S$, that is functors $F: \SmS^{\op} \rarr \Spc$, will be that of functors $\cF: \dStkFop \rarr \Top$. Since we work with Segal categories, we are therefore looking at the Segal category $\RHomdStkFopTop $ of all such functors, the categorical manifestation in $\Top = L \SetD$ of $\dStkF$. This is a Segal category from which the category of derived motivic spectra can be constructed. Using strictification (\cite{SeT}), which states that $\RHomdStkFopTop \cong L \Fun(\dStkFop, \SetD)$, we model objects \newline $\cF \in \RHomdStkFopTop$ by simplicial functors $\cF: \dStkFop \rarr \SetD$, and it is with those objects that we will be working. Here $L$ refers to Dwyer-Kan's simplicial localization functor (\cite{DK}), which produces a simplicial category, which we regard as a Segal category.\\ 

Recall how this is defined: if $B \in \CatD$, $V \subset B$ a subcategory, then one defines $L(B,V) = \diag F_*B [F_* V^{-1}] = (F_* B_n [ F_* V_n^{-1}])_n$, where the standard resolution $F_*C \in \CatD$ for some $C \in \Cat$, is defined by $F_nC = F^{n+1} C$, where $FC$ is the free category on $C$, with a generator $Fc$ for every non-identity morphism $c$ of $C$. This definition is important for us since it shows that objects of $LB$ and $B$ are the same and $V$-morphisms are sequentially crunched but other morphisms are otherwise left alone. We will use this quite a bit in projecting down notions of morphisms at the level of $\DkAfftet$ to the level of $\dSt_k = L \DkAfftet$, by essentially saying morphisms in $\dSt_k$ have a property $P$ if their lifts in $\DkAfftet$ do have this property. Here we use the shorthand $L\cC$ for $L(\cC,W)$.\\

\subsection{The category of derived motivic spectra}
We define the category of derived motivic spectra to be the category of motivic $\bP^1_X$-spectra in $\RHomdStkFopTop$ ($\PoneX$ will be defined later in subsection \ref{GenofP}). Objects $\cF$ therein are said to be motivic if they are $\bA^1$-local (or rather a generalization thereof as we will soon see) and are Nisnevich local in a sense defined below. The difference from the classical case is that we do not insist on our functors to satisfy Nisnevich descent, since for one thing this is classically equivalent to satisfying Nisnevich excision, but also because homotopy functors in Goodwillie calculus (\cite{GII}, \cite{GIII}) already satisfy a form of excision, and objects $\cF \in \RHomdStkFopTop$ being homotopy functors, it is therefore natural to ask that they satisfy a form of excision.\\

To define $\bP^1_X$, one needs a notion of $\bA^1$-localization as well as a notion of Nisnevich localization in our context, so we first define the generalization of $\bA^1$ in our setting, as well as the notion of Nisnevich square as well, along with the accompanying left Bousfield localizations.\\

\subsubsection{$\bA^1$-localization}
\begin{itemize}
	\item \underline{Generalization of $\bA^1$}.
		Recall from \cite{K} that a space $F$ over $S$ is $\bA^1$-homotopy invariant if for all smooth $X \in \SmS$, the canonical morphism of spaces $\Gamma(X,F) \rarr \Gamma(X \times \bA^1,F)$ is an isomorphism, where $\bA^1 = \Spec \bZ[t]$ is the usual affine line. Presently we generalize this to $\bA^1_k = \Spec k[x]$, which we regard as a constant simplicial object of $\skAlg^{\op}$. Then we let $\uAonek = \Ruh_{\bA^1_k}$, $\Ruh: \Ho(\skAlg^{\op}) \rarr \Ho(\DkAfftet)$ the Yoneda embedding, $\DkAfftet \subset \Fun(\skAlg, \SetD)$ the model category of stacks on the site $(\skAlg, \acute{e}t)$ (\cite{HAGII}).\\ 
	
	\item \underline{$\bA^1_k$-homotopy invariance}.
		We say $\cF : \dStkFop \rarr \Top$ is $\uAonek$-homotopy invariant if for all $X \in \dStkF$ its lift $\cF: \dStkFop \rarr \SetD$ satisfies $\cF(X) \simeq \cF(X \times \uAonek)$. We regard $\uAonek$ as an object of $\dSt_k$ by using the fact that $\Ho(\DkAfftet)$ and $\DkAfftet$ have the same objects, and the fact that $\dSt_k = L \DkAfftet$ and $\DkAfftet$ have the same objects as well.\\

	\item \underline{$\bA^1_k$-localization}.
		Going over the work done in \cite{K}, it is clear that one can recast it in an equivalent formalism using left Bousfield localizations after a bit of an enhancement. By \cite{Hi}, $\SetD$ is a left proper cellular model category, hence so is $\Fun(\dStkFop, \SetD):= \cM$, so $\cC$ being a class of morphisms therein, one can define the left Bousfield localization of this functor category with respect to $\cC$; it is a model structure on $\cM$ for which the weak equivalences are the $\cC$-local equivalences, cofibrations are those of $\cM$, and fibrations are those maps that have the right lifting property with respect to $\cC$-local equivalences that are also cofibrations. Also because $\cM$ is a model category, one can define a homotopy function complex $\map$ on it. Recall that if $\cM$ is a model category and $\cC$ is a class of morphisms in $\cM$, $W \in \Ob(\cM)$ is said to be $\cC$-local if it is fibrant and for all $f: A \rarr B$ in $\cC$, the induced map $\map(B,W) \rarr \map(A,W)$ is a weak equivalence in $\SetD$. Then a morphism $g:X \rarr Y$ in $\cM$ is said to be a $\cC$-local equivalence if for all $\cC$-local object $W$ of $\cM$, the induced map $\map(Y,W) \rarr \map(X,W)$ is a weak equivalence. Finally by Thm 4.1.1 of \cite{Hi}, fibrant objects of the left Bousfield localization $L_{\cC}\cM$ of $\cM$ with respect to $\cC$ are $\cC$-local objects, and objects of model categories have fibrant replacements. Thus $L_{\cC}\cM$ is effectively a model category of $\cC$-local objects, and $\Ho(L_{\cC}\cM)$ is consequently a category of $\cC$-local objects.\\

		It is easy to define $\uAonek$-local equivalences; if $\cC = \{ X \times \uAonek \rarr X \}$ is a class of morphisms in $\dStkF$, it gives rise to a class $h_{\cC} = \{ h_{X \times \uAonek} \rarr h_X \}$ in  $\RHomdStkFopTop$. Indeed, $\dStkF$ being a simplicial model category, it has its own simplicial mapping space $\Map_F : = \Map_{\dStkF}$. This allows us to define $h_X = \Map_F(-,X): \dStkFop \rarr \SetD$. As a first approximation an object $\cF$ of $\RHomdStkFopTop$ is said to be $\uAonek$-local if it is $\uAonek$-homotopy invariant, that is for all $X \in \dStkF$, $\cF(X) \rarr \cF(X \times \uAonek)$ is a weak equivalence of simplicial sets. Using Yoneda and the simplicial mapping space of $\Fun(\dStkFop, \SetD)$, this reads $\Map(h_X,\cF) \xrarr{\simeq} \Map(h_{X \times \uAonek}, \cF)$, which we generalize to the statement  $\map(h_X,\cF) \xrarr{\simeq} \map(h_{X \times \uAonek}, \cF)$, where map is the homotopy function complex of $\Fun(\dStkFop, \SetD)$. We take this as our definition of being $\uAonek$-local. Then a morphism $\cF \rarr \cF'$ is said to be an $\uAonek$-local equivalence if for all $\uAonek$-local objects $\cG$, $\map(\cF', \cG) \rarr \map(\cF,\cG)$ is a weak equivalence.\\
\end{itemize}

\subsubsection{Nisnevich localization}
\begin{itemize}
	\item \underline{Choice of Nisnevich excision over Nisnevich descent -  Goodwillie calculus}. We generalize the work done by Hoyois in \cite{H} and by Hall and Rydh in \cite{HR} on algebraic stacks to the level of derived stacks. The problem we encounter is that for a presheaf $\cF$ on stacks to satisfy Nisnevich descent, one naturally has to use Nisnevich coverings, which use the notion of residual gerbes that seem very difficult to generalize to functors on derived stacks. However, Hoyois shows that a presheaf $\cF$ on stacks satisfies Nisnevich descent if and only if it satisfies Nisnevich excision. Thus having this equivalence at the level of presheaves on algebraic stacks, we focus on functors $\cF: \dStkF^{\op} \rarr \Top$ that satisfy Nisnevich excision, in a sense to be precised.\\

		Presently we make an allusion to Goodwillie calculus (\cite{GII}, \cite{GIII}). Goodwillie considers homotopy functors, functors from spaces to spaces, or spectra, that preserve weak equivalences, and this would correspond to working with prestacks for us.  Such functors are said to be excisive if they map homotopy pushout squares to homotopy pullback squares, and this again is reminiscent of the Nisnevich excision condition. But excisive functors are really helpful in developing Calculus on functors, a functorial differential approach to working with prestacks, which is pertinent for us for the sake of deformations. In particular in \cite{GII}, Goodwillie shows that to any given homotopy functor $F$ one can associate what he calls a $n$-excisive functor $P_n F$ together with a map $F \rarr P_n F$, a sort of $n$-th Taylor polynomial for $F$, and those agree to order $n$. Further we have a tower of such polynomial approximations with a map $F \rarr \lim P_n F$. If $F$ is $\rho$-analytic, a condition based on excision, if $X \rarr Y$ is connected in a certain sense, then the connectivity of the map $F(X) \rarr (P_nF)(X)$ tends to infinity with $n$. Thus excision is at the heart of Goodwillie calculus. Since we work with motivic objects which are inherently differential in character, what we have here is yet another reason for choosing excision as a requirement for our functors.\\

We will adopt the definition of Nisnevich square of stacks by adapting the one provided in \cite{H} for algebraic stacks, that is we consider a cartesian square of stacks of the form:
		\beq
		\xymatrix{
			W \ar[d] \ar[r] &V \ar[d]^f \\
		U \ar[r]_e & X
			}\nonumber
		\eeq
		with $f$ \'etale, $e$ an open immersion with reduced complement $Z$ such that $Z \times_{X} V \xrarr{\cong} Z$. It turns out we will use Zariski open immersions for derived stacks, but otherwise all the other notions of morphisms will be preserved. For Hoyois, a presheaf of spaces, or spectra, $\cF$ on the category of algebraic stacks satisfies Nisnevich excision if the induced square:
		\beq
		\xymatrix{
			\cF(X) \ar[d]^{f^*} \ar[r]^{e^*} & \cF(U) \ar[d] \\
			\cF(V) \ar[r] & \cF(W)
			}\nonumber
		\eeq
		is homotopy cartesian. We will use the same formalism for functors $\cF: \dStkFop \rarr \Top$, by simply replacing open immersions by Zariski open immersions as defined in \cite{HAGII}. As usual we will work with representatives $\cF: \dStkFop \rarr \SetD$ of such functors by way of the strictification theorem  $\RHomdStkFopTop \cong L \Fun(\dStkFop, \SetD)$.
	\item \underline{Nisnevich  squares}. 		
		\begin{itemize}
			\item[-] \underline{Smooth morphisms}. Classically one works with Nisnevich squares in $\SmS$, which have $\dStkF$ as analog for us, so we seek to define a notion of smoothness in $\dSt_k$. First, the notions of smoothness in \cite{K} and \cite{HAGII} at the level of simplicial algebras coincide. Then one uses the model Yoneda embedding. By \cite{HAGII}, for $Q$ a class of morphisms in $(\skAlg)^{\op} = \DkAff$ stable by equivalences, homotopy pullbacks and compositions, the Yoneda embedding $\Ruh: \Ho( \DkAff) \rarr \Ho(\DkAfftet)$ allows one to extend the notion of $Q$-morphism to the level of representable stacks. $Q = \Sm$ is one such class by Def 1.2.7.2 of \cite{HAGII}, and by virtue of Prop 1.2.3.3 and Prop 1.2.7.3. Thus one can talk about smooth morphisms between representable stacks. Finally one uses the fact that $\dSt_k  = L (\DkAfftet)$. Thus smooth morphisms are well-defined in $\dStkF$.\\

			\item[-] \underline{Zariski open immersions}. In \cite{HAGII} it is argued that Zariski open immersions reduce to the classical notion of open immersion in the event that we put the trivial model structure on the category of commutative monoid objects in our base model categories. Thus Zariski open immersions appear to be the correct generalization in homotopical algebra of the classical notion of open immersion, and this is the notion we will use. By Prop 2.2.2.5 of \cite{HAGII} a morphism $f: A \rarr B$ in $\skAlg$ is a Zariski open immersion if the induced morphism of affine schemes $\Spec \pi_0 B \rarr \Spec \pi_0 A$ is a Zariski open immersion and $A \rarr B$ is strong, that is $\pi_* A \oT_{\pi_0 B} \pi_0 A \rarr \pi_* B$ is an isomorphism. Here we use the standard definitions from \cite{HAGII}: $|A| = \Map_{\skAlg}(1,A) \in \Ho(\SetD)$, $\pi_i A = \pi_i(|A|,*)$. To say that a morphism of affine schemes $\Spec B \rarr \Spec A$ is a Zariski open immersion means $A \rarr B$ is a Zariski open immersion in $\kAlg$, meaning it is flat, of finite presentation, and is such that $B \cong B \oT_A B$. Then we use the fact that we identify $\skAlg$ with the category of representable stacks on it to transfer that notion of Zariski open immersion at the level of $\skAlg$ to the level of stacks. To be precise, if $Q$ is a class of morphisms in $\skAlg^{\op}$ stable by equivalences, homotopy pullbacks and compositions, then the Yoneda embedding $\Ruh: \Ho(\skAlg^{\op}) \rarr \Ho(\DkAfftet)$ allows one to lift the notion of morphism in $Q$ at the level of $\skAlg$ to the level of morphism of representable stacks, a fact we have used before. This is true of Zariski open immersions, which by Def 1.2.6.7 of \cite{HAGII} are finitely presented and are formal Zariski open immersions. Both classes are stable by equivalences, compositions and homotopy pushouts by Prop 1.2.3.3 and Prop 1.2.6.3. Thus we can talk of Zariski open immersions between representable stacks. We use the notation $\RuSpec A = \Ruh_{\Spec A}$ for $A \in \skAlg$. \\

				Observe that Toen and Vezzosi define a morphism of $D^{-}$-stacks to be a Zariski open immersion if it is a locally finitely presented and flat monomorphism. Since these notions do not allow us to find a transparent notion of reduced complement for the sake of defining a Nisnevich square of stacks, we use representable stacks instead for simplicity. If this seems to be a bit restrictive, recall that motives can be defined starting from schemes, which through their functors of points correspond to representable stacks, thus working with representable stacks within $\dStkF$ is sufficient for our purposes.\\

				\newpage

			\item[-] \underline{reduced complement}. For the sake of having a Nisnevich square, we need the notion of a reduced complement of a Zariski open immersion of representable stacks. Following \cite{TV5}, if $A \rarr B$ is a finitely presented morphism of commutative rings, then $\Spec B \rarr \Spec A$ is an open immersion if and only if the restriction functor $D^{-}(B) \rarr D^{-}(A)$ between derived categories is fully faithful. Now starting from $\RuSpec B \rarr \RuSpec A$ being a Zariski open immersion of representable stacks, means $\Spec B \rarr \Spec A$ is a Zariski open immersion, i.e. $A \rarr B$ is a Zariski open immersion in $\skAlg$, so it is strong and $\Spec \pi_0 B \rarr \Spec \pi_0 A$ is a Zariski open immersion, that is $ \pi_0 A \rarr \pi_0 B$ is a Zariski open immersion, so finitely presented and a formal Zariski open immersion by Def 1.2.6.7 of \cite{HAGII}, which means flat and $f_*: \Ho(\pi_0B \text{-Mod}) \rarr \Ho(\pi_0 A \text{-Mod})$ fully faithful, inducing a fully faithful map at the level of derived categories. It follows $\Spec \pi_0 B \rarr \Spec \pi_0 A$ is an open immersion.\\
				 
				To recap, $\RuSpec B \rarr \RuSpec A$ Zariski open immersion means $\Spec \pi_0 B \rarr \Spec \pi_0 A$ is an open immersion, so $| \Spec \pi_0 B| \cong U^{\text{open}} \subset \Spec \pi_0 A$, and it has a closed complement $Z = | \Spec \pi_0 A| - U$ on which we can put a reduced closed subscheme structure.\\

				Now because $Z$ is a closed subscheme of an affine scheme $\Spec \pi_0 A$, it is affine itself, $Z = \Spec C$, regard $C = \pi_0 C_{\bullet}$, $C_{\bullet} \in \skAlg$ being $C$ viewed as a constant simplicial object of $\kAlg$, with accompanying representable stack $\RuSpec C_{\bullet}$, along with a morphism to $\RuSpec A$, which we regard as the reduced complement to the Zariski open immersion $\RuSpec B \rarr \RuSpec A$.\\

			\item[-] \underline{\'etale morphisms}.
				By Def 1.2.6.7 of \cite{HAGII}, \'etale morphisms are finitely presented and formally \'etale, and by Prop 1.2.3.3 and Prop 1.2.6.3 respectively, those are stable by equivalences, compositions and homotopy pushouts, thus we can lift the notion of \'etale morphism at the level of $\skAlg$ to the level of representable stacks in $\DkAfftet$ by using the usual argument. Finally we invoke the fact that $\dSt_k = L \DkAfftet$ to argue that \'etale morphisms are well-defined in $\dStkF$.\\
		\end{itemize}
	\item \underline{Nisnevich localization}.
In our present situation of working with functors $\cG: \dStkFop \rarr \SetD$, in a first approximation $\cG$ is said to be Nis-local if it satisfies Nisnevich excision, and a morphism $\cF \rarr \cF'$ of objects of $\Fun(\dStkFop, \SetD)$ is said to be a Nis-local equivalence if for all Nis-local $\cG$ in this functor category, $\map(\cF',\cG) \rarr \map(\cF,\cG)$ is a weak equivalence of simplicial sets.\\

It is more difficult to define a left Bousfield localization with respect to Nisnevich squares than it was to define a left Bousfield localization with respect to maps of the form $X \times \uAonek \rarr X$. The difficulty resides in taking as our class of morphisms for left Bousfield localization purposes a class of squares. We show this can be done presently.\\

Recall that to satisfy Nisnevich excision means that for any Nisnevich square:
\beq
\xymatrix{
	W \ar[d] \ar[r] &V \ar[d] \\ 
	U \ar[r] & X
} \nonumber
\eeq
the square:
\beq
\xymatrix{
	\cF(X) \ar[d] \ar[r] & \cF(U) \ar[d] \\
	\cF(V) \ar[r] & \cF(W)
} \nonumber
\eeq
is homotopy cartesian, which means $\cF(X) \xrarr{\simeq} \cF(V) \times^h_{\cF(W)} \cF(U)$. Now Thm 18.1.10 of \cite{Hi} states that if $\cM \in \CatD$, $\cC$ is a small category, $\cX$ is a $\cC$-diagram in $\cM$, $Y$ is an object of $\cM$, then we have an isomorphism $\Map(\hocolim_{\cC} \cX,Y) \cong \holim_{\cCop}\Map(\cX,Y)$, $\Map$ being the simplicial mapping space of $\cM$. Applying this to $\cM = \Fun(\dStkFop, \SetD)$, with $\cC = \{ U \larr W \rarr V \}$, $Y = \cF$ and $h = \cX: \cC \rarr \cM$, we have $\cF(U) \times^h_{\cF(W)}\cF(V) \simeq \cF(V \coprod^h_W U)$. Thus we see that $\cF$ satisfying Nisnevich excision can be recast as saying $\cF$ is local with respect to maps in $\cC = \{ V \coprod^h_W U \rarr X \}$ induced by Nisnevich squares. Thus satisfying Nisnevich excision is equivalent to being local with respect to Nisnevich squares, which we just refer to as being Nis-local. This is a first approximation. For the sake of using left Bousfield localizations we generalize this as follows. Let $h_{\cC} = \{ h_{U \coprod^h_W V} \rarr h_X \}$ be the class of morphisms of $\Fun(\dStkFop, \SetD)$ induced by the above class $\cC = \{ U \coprod^h_W V \rarr X \}$ in $\dStkF$, and replace the statement:
\beq
\Map(h_X,\cF) \simeq \cF(X) \xrarr{\simeq} \cF(U \coprod^h_W V) \simeq \Map(h_{V \coprod^h_W U}, \cF) \nonumber
\eeq
by $ \map(h_X,\cF) \xrarr{\simeq} \map(h_{V \coprod^h_W U}, \cF) $. We take this as our definition of being Nis-local. Then $\cF \rarr \cF'$ is a Nis-local equivalence if for any Nis-local object $\cG$, we have $\map(\cF',\cG) \rarr \map(\cF,\cG)$ is a weak equivalence of simplicial sets.\\
\end{itemize}

\subsubsection{ Generalization of $\bP^1_S$} \label{GenofP}
Recall how $\bP^1_S$ was constructed in \cite{K}: $\bP^1_S = (M_S(\bP^1_S),\infty)$, where if $X \in \SmS$, $M_S(X) =\Lmot(h_S(X))$, where for a space $F$ over $S$, $\Lmot(F) = \lim(\LAone \circ \LNis)^n(F)$ where $\LAone$ is an accessible localization onto $\bA^1$-homotopy invariant spaces, $\LNis$ is the left exact localization onto Nisnevich local spaces. We also have $h_S(X) = \Map_S(-,X): \SmS^{\op} \rarr \Spc$. Going over the construction of $\bP^1_S$ in \cite{K}, it becomes clear that one can generalize it to our setting where we work with derived schemes in $\dSchk$, thus $\bP^1_{iX} \in \dSchk$ is well-defined.\\

We then use the fact that $\dSchk \in \CatD$, hence it has a simplicial mapping object, hence so does $(\dSchk)_{/iX}$. We denote its simplicial mapping space by $\Map_X$. This allows us to see $h_{\PoneiX} = \Map_X(-,\bP^1_{iX}): (\dSchk)_{/iX}^{\op} \rarr \SetD$ as an object of $\dStkF$. Now we need to see this as a functor $\dStkFop \rarr \SetD$. This is achieved by taking its functor of points $\Map_F(-,h_{\PoneiX}) = h_F \circ h_{\PoneiX}$. We then take the motivic localization $\Lmot = \LuAonekNis$ of this object. Here $\LuAonekNis$ is a left Bousfield localization of $\Fun(\dStkFop, \SetD)$ with respect to the class $\{ h_{X \times \uAonek} \rarr h_X \} \cup \{ h_{U \coprod^h_W V} \rarr h_X \}$, the first subclass being informally referred to as $\uAonek$, and the second one by $Nis$, and again it is induced by Nisnevich squares $W \rrarr U,V \rrarr X$. \\

We define the category of derived motivic spaces by:
\beq
\Mot(F) = \Ho(\LuAonekNis \SetD^{\dStkFop}) \nonumber
\eeq
We regard $\Lmot h_F \circ h_{\PoneiX}$ as an object of this category, and we still denote it by $\bP^1_X$ for simplicity.\\

\subsubsection{Localizations}
Presently we have the following picture:
\beq
\xymatrix{
	\SetD^{\dStkFop} \ar[r]^{\LuAonekNis} &\LuAonekNis \SetD^{\dStkFop} \ar[r]^-{\Ho} & \Mot(F) \ar[d]^{L} \\
	\RHomdStkFopTop \ar[u]^{L^{-1}} && \SeMot(F)
} \nonumber
\eeq
where $\SeMot(F) = L \Mot(F)$ is the Segal category of derived motivic spaces on $F$, whose objects are also those of $\RHomdStkFopTop$. That this holds is by virtue of the fact that the equivalences in $\Fun(\dStkFop, \SetD)$ used in the simplicial localization on the left in the above diagram are not involved in the left Bousfield localization $\LuAonekNis$, hence are preserved throughout, and can be used later as in the right hand side. This also means that objects of the Segal category $\SeMot(F)$ can by modeled in $\Mot(F)$ by functors $\cF: \dStkFop \rarr \SetD$. We will use this fact in the next subsection.\\

\subsubsection{$\cT$-spectra}
Recall from \cite{K} that objects of the category $\SHT(S)$ of motivic $\cT$-spectra are sequences $\{ F_n \}_{n \geq 0}$ of fibred pointed spaces with structural isomorphism $F_n \rarr \Omega_{\cT}(F_{n+1})$, where $\Omega_{\cT} = \uHom(\cT,-)$, using the internal hom of the category of spaces over $S$. We mention in passing the fact that we will not work with pointed spaces, since we will not need it in our work, but the conscientious reader will have no difficulty making the proper changes for the pointed setting.\\

We generalize this in our context by defining the category of derived motivic $\cT$-spectra, denoted $\MotSp_{\cT}(F)$, $F = h_{iX} \in \dSt_k$ if $X \in \Schk$, with $\cT$ an object of $\SeMot(F)$. The objects of $\MotSp_{\cT}(F)$ are sequences $\{ \cF_n \}_{n \geq 0}$ of objects of $\SeMot(F)$ along with structural isomorphisms $\cF_n \cong \Omega_{\cT}(\cF_{n+1})$, $\Omega_{\cT} = \uHom_F(\cT,-)$, $\uHom_F$ the internal Hom of $\SeMot(F)$, that we define as follows.\\

$\SetD$ is a symmetric monoidal model category with an internal hom defined as follows. For $K,L \in \SetD$, one defines the function complex $\Map(K,L)$ in degree $n$ by $\Map(K,L)_n = \Hom_{\SetD}(K \times \Delta_n, L)$. For $\cF, \cG, \cH$ three objects of $\SeMot(F)$ modeled as functors $\dStkFop \rarr \SetD$, we have:
\beq
\Hom_{\SetD}(\cF(U) \times \cG(U), \cH(U)) \cong \Hom_{\SetD}(\cF(U), \Map(\cG(U), \cH(U)) \nonumber
\eeq
for $U \in \dStkF$. For a given functor $\cF: \cCop \rarr \cD$, define $\cF^{/\op}: \cC \rarr \cD$. Now define $\cF \oT \cG$ and $\uHom(\cG,\cH)$ pointwise by $(\cF \oT \cG)(U) = \cF(U) \times \cG(U)$, and $\uHom(\cG,\cH)(U) = \Map(\cG^{/\op}(U),\cH(U))$. In this manner we have $\uHom(\cG,\cH)$ and $\cF \oT \cG$ being functors $ \dStkFop \rarr \SetD$. These being defined, we have:
\beq
\Hom_{\Mot(F)} ( \cF \oT \cG,\cH) \cong \Hom_{\Mot(F)}(\cF, \uHom(\cG,\cH)) \nonumber
\eeq
which projects down to the same statement at the level of $\SeMot(F)$ after simplicial localization, if we look at individual elements of these hom sets, and we use the fact that objects of $\Mot(F)$ are also those of $\SeMot(F)$. Thus:
\beq
\Hom_{\SeMot(F)}( \cF \oT \cG,\cH) \cong \Hom_{\SeMot(F)}(\cF, \uHom(\cG,\cH)) \nonumber
\eeq

\subsubsection{ Construction of $\MotSp(X)$.}
We define $\MotSp(X) = \MotSp_{\PoneX}(F)$ to be the category of derived motivic spectra on $X \in \Schk$, with $F = h_{iX}$, whose objects are $\PoneX$-spectra of functors $\cF: \dStkFop \rarr \SetD$ - modeling objects of $\RHomdStkFopTop$- that are $\uAonek$-local and $Nis$-local. $\MotSp(X)$ appears first as a delocalization, followed by a left Bousfield localization, upon which we take the corresponding homotopy category, followed by a simplicial localization, to finish with a spectralization, as in the diagram below:
\beq
\xymatrix{
	\SetD^{\dStkFop} \ar[r]^{\LuAonekNis} &\LuAonekNis \SetD^{\dStkFop} \ar[r]^-{\Ho} & \Mot(F) \ar[d]^{L} \\
	\RHomdStkFopTop \ar[u]^{L^{-1}} && \SeMot(F) \ar[d]^{\Sp_{\PoneX}} \\
	&& \MotSp(X)
} \nonumber
\eeq
Looking ahead, if one is interested in making contact with the theory of motives, one has to find an appropriate topology on $\dStkFop$ to define stacks in $\MotSp(X)$, take the derived category thereof, take the analog in our situation of a $\DAone$ localization, producing the equivalent for us of the derived categories of sheaves of \cite{CD}, and possibly use an appropriate reduction to get a category equivalent to $DM_B(X)$. Observe that this subsumes developing a 6-functors formalism in the present setting, something we have not touched upon.

\bigskip
\footnotesize
\noindent
\textit{e-mail address}: \texttt{rg.mathematics@gmail.com}.

\end{document}